\documentclass[11pt]{article}

\usepackage{amsfonts}
\usepackage{amsmath}
\usepackage{amssymb}
\usepackage{latexsym,color}
\usepackage{amsbsy}
\usepackage{graphicx}

\newtheorem{thm}{Theorem}[section]
\newtheorem{con}{Conjecture}[section]
\newtheorem{lemma}[thm]{Lemma}
\newtheorem{cor}[thm]{Corollary}
\newtheorem{pro}[thm]{Proposition}
\newtheorem{example}[thm]{Example}
\newtheorem{definition}[thm]{Definition}

\newtheorem{remark}[thm]{Remark}
\newtheorem{Algorithm}[thm]{Algorithm}

\newcommand{\comment}[1]{}

\newcommand{\N}{{\mathbb N}}
\newcommand{\C}{{\mathbb C}}

\newcommand{\pr}{\parallel}
\newcommand{\spec}{\mbox{spec}}

\newcommand{\ncom}{\newcommand}
\ncom{\ns}{\normalsize}
\ncom{\la}{\lambda}
\ncom{\bm}{\boldmath}
\ncom{\noi}{\noindent}
\ncom{\bq}{\begin{equation}}
\ncom{\eq}{\end{equation}}
\ncom{\beqn}{\begin{eqnarray*}}
\ncom{\eeqn}{\end{eqnarray*}}
\ncom{\ba}{\begin{array}}

\ncom{\ea}{\end{array}}
\ncom{\beq}{\begin{eqnarray}}
\ncom{\eeq}{\end{eqnarray}}
\ncom{\nno}{\nonumber}
\ncom{\hs}{\mbox{\hspace{.25cm}}}
\ncom{\rar}{\rightarrow}
\ncom{\Rar}{\Rightarrow}
\ncom{\noin}{\noindent}
\ncom{\bc}{\begin{center}}
\ncom{\ec}{\end{center}}
\ncom{\sz}{\scriptsize}
\ncom{\fpd}{\Phi(\pi^{'})}
\ncom{\fp}{\Phi(\pi) }
\ncom{\nk}{\left< \begin{array}{c}
                       n\\k \end{array} \right>}
\ncom{\nd}{1^{'},2^{'},\cdots,n^{'}}
\ncom{\R}{I\!\!R}
\ncom{\de}{\bigtriangleup (F_{2n},\leq)}
\ncom{\del}{\bigtriangleup}
\ncom{\cov}{<\!\!\!\!\cdot }
\ncom{\bt}{\begin{thm}}
\ncom{\bcon}{\begin{con}}
\ncom{\et}{\end{thm}}
\ncom{\econ}{\end{con}}
\ncom{\bl}{\begin{lemma}}
\ncom{\el}{\end{lemma}}
\ncom{\bco}{\begin{cor}}
\ncom{\ds}{\displaystyle}
\ncom{\eco}{\end{cor}}
\ncom{\bp}{\begin{pro}}
\ncom{\ep}{\end{pro}}
\ncom{\bex}{\begin{example}}
\ncom{\eex}{\end{example}}
\ncom{\bd}{\begin{definition}}
\ncom{\ed}{\end{definition}}
\ncom{\brm}{\begin{remark}}
\ncom{\erm}{\end{remark}}
\ncom{\bal}{\begin{Algorithm}}
\ncom{\eal}{\end{Algorithm}}
\ncom{\ol}{\overline}

\ncom{\pf}{\noi {\bf Proof  }}
\ncom{\be}{\begin{enumerate}}
\ncom{\ee}{\end{enumerate}}
\ncom{\s}{\subset}
\ncom{\T}{{\cal T}}
\ncom{\B}{{\cal B}}
\ncom{\A}{{\cal A}}

\title{\Large{\textcolor{black} {\bf Symmetric chains, Gelfand-Tsetlin
chains, and the Terwilliger algebra of the binary Hamming scheme}}}

\author{{\textcolor{black} {\bf Murali K. Srinivasan}} \\
{\em  \normalsize{Department of Mathematics}}\\
{\em  \normalsize{Indian Institute of Technology, Bombay}}\\
{\em  \normalsize{Powai, Mumbai 400076, INDIA}}\\
{\bf  \texttt{mks@math.iitb.ac.in}}}

\begin{document}
\date{}
\maketitle
\begin{center}{\em To the memory of Mobi}\end{center}

\begin{abstract}
The de Bruijn-Tengbergen-Kruyswijk (BTK) construction is a
simple algorithm that produces an explicit symmetric chain
decomposition of a product of chains. We linearize the BTK algorithm and
show that it produces an explicit symmetric Jordan basis (SJB).
In the special case of a Boolean algebra
the resulting SJB is orthogonal with respect to the standard inner
product and, moreover, we can write down an explicit formula for the ratio
of the lengths of the successive vectors in these chains (i.e., the singular
values). This yields a new, constructive proof of the 
explicit block diagonalization of the Terwilliger algebra
of the binary Hamming scheme.  
We also give a representation theoretic characterization of this basis that
explains its orthogonality, namely, that it is the canonically defined (upto
scalars)
symmetric Gelfand-Tsetlin basis.
\end{abstract}

{\textcolor{black}{ \section{\bf  Introduction }}}  

The de Bruijn-Tengbergen-Kruyswijk (BTK) construction  is a simple visual
algorithm in matching theory that produces an explicit symmetric chain decomposition of a (finite)
product of (finite) linear orders. We show that the BTK algorithm admits a
simple and natural linear analog. The main purpose of this paper is to
study the linear BTK algorithm as an object in itself. It  
enables us to explicitly study the up operator (i.e., the nilpotent operator
taking an element to the sum of the elements covering it)
on a product of linear orders by producing an explicit symmetric 
Jordan basis (SJB). In the special case of a Boolean algebra
the resulting SJB is orthogonal with respect to (wrt) the standard inner
product and, moreover, we can write down an explicit formula for the ratio
of the lengths of the successive vectors in these chains (i.e., the singular
values). This yields a new constructive proof
of the explicit block diagonalization of the Terwilliger algebra
of the binary Hamming scheme, recently achieved by Schrijver.
We also give a representation theoretic characterization of this basis that
explains its orthogonality, namely, that it is the canonically defined 
(upto scalars)
symmetric Gelfand-Tsetlin basis (wrt the up operator on 
the Boolean algebra).

A (finite) {\em graded poset} is a (finite) poset $P$ together with a 
{\em rank function} 
$r: P\rar \N$ such that if $q$ covers $p$ in $P$ then 
$r(q)=r(p)+1$. The {\em rank} of $P$ is $r(P)=\mbox{max }\{r(p): p\in P\}$
and,
for $i=0,1,\ldots ,r(P)$, $P_i$ denotes the set of elements of $P$ of rank
$i$. A {\em symmetric chain} in a graded poset $P$ is a sequence
$(p_1,\ldots ,p_h)$ of elements of $P$ such that $p_i$ covers $p_{i-1}$, for
$i=2,\ldots h$, and 
$r(p_1) + r(p_h) = r(P)$, if $h\geq 2$, or else $2r(p_1) = r(P)$, if $h=1$.  
A {\em symmetric chain
decomposition} (SCD) of a graded poset $P$ is a decomposition of $P$ into pairwise disjoint
symmetric chains. 

We now define the linear analog of a SCD. 
For a finite set $S$, let $V(S)$ denote the complex vector space with $S$ as
basis. Let $P$ be a graded poset with $n=r(P)$. Then we have 
$ V(P)=V(P_0)\oplus V(P_1) \oplus \cdots \oplus V(P_n)$ (vector space direct
sum).
An element $v\in V(P)$ is {\em homogeneous} if $v\in V(P_i)$ for some $i$,
and we extend the notion of rank to homogeneous elements by writing $r(v)=i$.
A linear map $T:V(P)\rar V(P)$ is said to be {\em order raising} if, for all
$p\in P$, $T(p)$ is a linear combination of the elements covering $p$ (note
that this implies that $T(p)=0$ for all maximal elements of $P$ and that
$T(v)$ is homogeneous for homogeneous $v$). The {\em
up operator}  $U:V(P)\rar V(P)$ is defined, for $p\in P$, by
$U(p)= \sum_{q} q$,  
where the sum is over all $q$ covering $p$.
Let $T$ be an
order raising map on a graded poset $P$. 
A {\em graded Jordan chain} in $V(P)$ with respect to $T$ (wrt $T$ for
short) is a sequence
$v=(v_1,\ldots ,v_h)$ of nonzero homogeneous elements of $V(P)$ 
such that $T(v_{i-1})=v_i$, for
$i=2,\ldots h$, and $T(v_h)=0$ (note that the
elements of this sequence are linearly independent, being nonzero and of
different ranks). We say that $v$ {\em 
starts} at rank $r(v_1)$ and {\em ends} at rank $r(v_h)$. 
If, in addition, $v$ is symmetric, i.e., $r(v_1) + r(v_h) = r(P)$, if $h\geq
2$, or else $2r(v_1)= r(P)$, if $h=1$, we say that $v$ is a  
{\em symmetric Jordan chain}.  
A {\em graded Jordan basis}  of $V(P)$ wrt $T$ is a basis of $V(P)$
consisting of a disjoint union of graded Jordan chains in $V(P)$ wrt $T$.
If every chain in a graded Jordan basis is symmetric we speak of a {\em
symmetric Jordan basis} (SJB) of $V(P)$ wrt $T$. 
When $T$ is the up map $U$, we drop the ``wrt $U$" from the notation
and speak of a graded Jordan basis/SJB of $V(P)$.

Let $n$ be a positive integer and let $k_1,\ldots ,k_n$ be nonnegative 
integers. Define
$$M(n,k_1,\ldots ,k_n) = \{(x_1,\ldots ,x_n)\in \N^n\;:\;0\leq x_i \leq k_i,\mbox{
for all }i\},$$
and partially order it by componentwise $\leq$. The cardinality of
$M(n,k_1,\ldots ,k_n)$
is $(k_1+1)\cdots (k_n+1)$ and the rank of $(x_1,\ldots , x_n)$ is 
$x_1+\cdots +x_n$, so $r(M(n,k_1,\ldots ,k_n))=k_1 + \cdots + k_n$.
It is easily seen that $M(n,k_1,\ldots ,k_n)$ is (order) isomorphic to a product
of $n$ chains of lengths $k_1,\ldots ,k_n$, respectively. Two special cases of 
$M(n,k_1,\ldots ,k_n)$ are of interest: the {\em uniform} case, where
$k_1=\cdots =k_n=k$ and we write $M(n,k)$ for $M(n,k,\ldots ,k)$ and the
{\em Boolean algebra or set} case, where $k_1=\cdots=k_n=1$ 
and we write $B(n)$ for $M(n,1)$.   

An algorithm to
construct an explicit SCD of $M(n,k_1,\ldots ,k_n)$ was given by
de Bruijn, Tengbergen, and Kruyswijk {\bf\cite{btk,a,e}}. We call this the BTK
algorithm. 
Canfield {\bf\cite{c}}, Proctor {\bf\cite{p}}, and Proctor, Saks, and
Sturtevant {\bf\cite{pss}} proved the existence of a SJB of
$V(M(n,k_1,\ldots ,k_n))$. These authors work in the more general context of
Sperner theory which we do not need here. An overview of this area is given
in Chapter 6 of 
Engel's book {\bf\cite{e}}. 
In Section 3 we present a linear analog of the  BTK algorithm 
and prove the following result.

\bt \label{mt1} For a positive integer $n$ and nonnegative integers
$k_1,\ldots ,k_n$, the linear BTK algorithm
constructs an explicit SJB of $V(M(n,k_1,\ldots ,k_n))$. The vectors in this
basis have integral coefficients when expressed in the standard basis
$M(n,k_1,\ldots ,k_n)$.
\et


When applied to the Boolean algebra $B(n)$ the linear BTK algorithm
has properties that go well  
beyond producing an SJB of $V(B(n))$. 
We now discuss this.
Perhaps the linear BTK algorithm (or a variant) 
has interesting properties also  in the 
general case but we are unable to say anything on this point here.

When the linear BTK algorithm is run  on $B(n)$, the 
resulting SJB's are all
orthogonal wrt the standard inner product on $B(n)$. Moreover, any two
symmetric
Jordan chains starting at rank $k$ and ending at rank $n-k$ ``look alike" in the sense
made precise in the following result.  
Let $\langle , \rangle$ denote the standard inner product on $V(B(n))$, i.e., 
$\langle X,Y \rangle =
\delta (X,Y)$, (Kronecker delta) for $X,Y\in B(n)$.
The {\em length} $\sqrt{\langle v, v \rangle }$ of $v\in V(B(n))$ is denoted 
$\pr v \pr$.
The following result is proved in Section 3. (In the formulation below item
(ii) is clearly implied by item (iii) but for convenience of 
later reference we have spelt out item (ii) explicitly.)

\bt \label{mt2}
Let $O(n)$ be the SJB produced by the linear BTK algorithm 
when applied to $B(n)$.

\noi (i) The elements of $O(n)$ are orthogonal with respect to 
$\langle , \rangle$.

\noi (ii) Let $0\leq k \leq \lfloor n/2 \rfloor$ and let 
$(x_k,\ldots ,x_{n-k})$  and $(y_k,\ldots ,y_{n-k})$ be any two symmetric Jordan chains
in $O(n)$ starting at rank $k$ and ending at rank $n-k$. Then
$$\frac{\pr x_{u+1} \pr}{\pr x_u \pr} =
\frac{\pr y_{u+1} \pr}{\pr y_u \pr},\;\;k\leq u < n-k.$$

\noi (iii) In the notation of part (ii) we have, for $k\leq u < n-k$, 
\beq \label{mks}
\frac{\pr x_{u+1} \pr}{\pr x_u \pr} & = &
\sqrt{(u+1-k)(n-k-u)}\\
\label{sch}
  & = & (n-k-u){{n-2k}\choose{u-k}}^{\frac{1}{2}}
                {{n-2k}\choose{u+1-k}}^{-\frac{1}{2}}.
\eeq
\et
In a recent
breakthrough, Schrijver {\bf\cite{s}} obtained new polynomial time
computable upper bounds on binary
code size using semidefinite programming and the Terwilliger algebra (this
approach was later extended to nonbinary codes in {\bf\cite{gst}}). There
are two main steps involved here. The first is to (upper) bound 
binary code size by the optimal value of an exponential size semidefinite
program and the second is to reduce the semidefinite program to polynomial
size by explicitly block diagonalizing the Terwilliger algebra.
For
background on coding theory we refer to {\bf\cite{gst,s}}. In this paper we
consider the second step. Theorem \ref{mt2} contains most of the information
necessary to explicitly
block diagonalize the Terwilliger algebra of the binary Hamming scheme.
We show this in Section 2.
Here we would like to add a few remarks
about the present proof of explicit block diagonalization. There are two
proofs available for this result: the 
linear algebraic proof of Schrijver and the representation theoretic 
proof of Vallentin 
{\bf\cite{v}} based on the work of 
Dunkl {\bf\cite{d1,d2}}. Our proof can be seen as a constructive version of Schrijver's proof
that also has representation theoretic meaning (see Theorem \ref{mt3}
below).
The basic pattern
of the proof is the same as in {\bf\cite{s}}. Given Theorem \ref{mt2}, the
rest of the proof  is a binomial inversion argument from
{\bf\cite{s}}. On the other hand, though not explicitly stated in this form
in {\bf\cite{s}}, the existence of a SJB of $V(B(n))$ 
satisfying (i), (ii), and
(iii) of Theorem \ref{mt2} easily follows from the results in {\bf\cite{s}}.
So the new ingredient here is the explicit construction of the SJB $O(n)$ and
its representation theoretic characterization in Theorem \ref{mt3} below. We
remark here that this explicit construction is primarily of mathematical
interest and is not important from the complexity point of view since even
to write down $O(n)$ takes exponential time.
We 
rewrite equation  (\ref{mks}) as equation (\ref{sch}) so that the final
formula for the  block diagonalization turns out to equal Schrijver's which
is in a very convenient form with respect to the location of square roots 
(see Theorem \ref{scht}).
 
Our proof of Theorem \ref{mt2} is self-contained and elementary but more insight into the
result is obtained by using a bit of representation theory. We first give a short 
proof of
existence of a SJB of $V(B(n))$ satisfying (i), (ii), and (iii) of Theorem
\ref{mt2} by putting together standard and well-known facts from representation theory.

Define the down operator $D$ on $V(B(n))$ analogous to the up operator and
define the operator $H$ on $V(B(n))$ by $H(v_i) = (2i -n)v_i,$ $v_i\in
V(B(n)_i),\;i=0,1,\ldots, n$. It is easy to check that $[H,U]=2U$,
$[H,D]=-2D$, and $[U,D]=H$. Thus the linear map ${\mbox{sl}}(2,\C)\rar
{\mbox{gl}}(V(B(n)))$ given
by 
$$\left(\ba{cc} 0&1\\0&0\ea \right) \mapsto U,\;
\left(\ba{cc} 0&0\\1&0\ea \right) \mapsto D,\;
\left(\ba{cc} 1&0\\0&-1\ea \right) \mapsto H$$
is a representation of ${\mbox{sl}}(2,\C)$. Decompose $V(B(n))$ into irreducible
${\mbox{sl}}(2,\C)$-submodules and let
$W$ be an irreducible in this decomposition
with dimension $l+1$. It follows from the representation theory of
${\mbox{sl}}(2,\C)$ (see Section 2.3.1 in {\bf\cite{gw}}) that there exists a basis
$\{v_0,v_1,\ldots ,v_l\}$ of $W$ such that, for $i=0,1,\ldots ,l$, we have (below
we take $v_{-1}=v_{l+1}=0$)
\beq
\label{sl2c} U(v_{i})=v_{i+1},\;
D(v_{i})=i(l-i+1)v_{i-1},\;
H(v_{i})=(2i-l)v_{i}&&
\eeq
So the eigenvalues of $H$ on $v_0,v_1,\ldots ,v_l$ are, respectively,
$-l,-l+2,\ldots ,l-2,l$.
It now follows from the definition of $H$ that each $v_i$ is homogeneous and 
that $(v_0,\ldots ,v_l)$ is a
symmetric Jordan chain in $V(B(n))$. Thus there exists a SJB of $V(B(n))$.
Note also that the basis $\{v_0,\ldots ,v_l\}$ of $W$ is canonically
determined (upto a common scalar multiple) since the eigenvalues of the
$v_i$ on $H$ are distinct. Let $r(v_0)=k$ and put $x_j=v_{j-k},\;k\leq j
\leq n-k$. The symmetric Jordan chain $(v_0,\ldots ,v_l)$ now gets rewritten
as $(x_k,\ldots ,x_{n-k})$ and (\ref{sl2c}) becomes, for $k\leq u \leq n-k$,
\beq
\label{sl2c1} U(x_u)=x_{u+1},\;
D(x_u)=(u-k)(n-k-u+1)x_{u-1}.
\eeq
Now we use the substitution action of the  symmetric group $S_n$ on $B(n)$.
As is easily seen the
existence of a SJB of $V(B(n))$ satisfying (i) and (ii) of Theorem \ref{mt2}
follows from the following facts by an application of Schur's lemma:

\noi (a) Existence of some SJB of $V(B(n))$.

\noi (b) $U$ is $S_n$-linear.

\noi (c) For $0\leq k \leq n$, $V(B(n)_k)$ is the sum of
$\mbox{min }\{k,n-k\}+1$ distinct
irreducible $S_n$-modules (this result is well known).

\noi (d) For a finite group $G$, a $G$-invariant inner product on an irreducible $G$-module
 is unique upto scalars.

We now   again use the ${\mbox{sl}}(2,\C)$ action to prove  Theorem \ref{mt2}(iii).
Let $J(n)$ be a SJB of $V(B(n))$ satisfying parts (i) and (ii) of Theorem
\ref{mt2}. Normalize $J(n)$ to get an
orthonormal basis $J'(n)$ of $V(B(n))$. Let $(x_k,\ldots ,x_{n-k})$, $x_u
\in V(B(n)_u)) \mbox{ for all } u$, be a symmetric Jordan chain in $J(n)$. Put
$x_u' = \frac{x_u}{\pr x_u \pr}$ and $\alpha_u = \frac{\pr x_{u+1} \pr}{\pr
x_{u} \pr},\;k\leq u \leq n-k$ (we take $x_{k-1}=x_{n-k+1}=0$). We have,
for $k\leq u \leq n-k$,
\beq \label{trick}
U(x_{u}')=\frac{U(x_{u})}{\pr x_{u} \pr}=\frac{x_{u+1}}{\pr x_{u}
\pr}=\alpha_{u} x_{u+1}'.&&
\eeq
Now observe that the matrices, in the standard basis, 
of $U$ and $D$ are real and transposes of
each other. Since $J'(n)$ is orthonormal
wrt the standard inner product, it follows that the matrices of
$U$ and $D$, in the basis $J'(n)$, must be adjoints of each other.
Thus  we must have, using (\ref{trick}),
$D(x_{u+1}')=\alpha_u x_{u}'$. So the subspace spanned by $\{x_k,\ldots
,x_{n-k}\}$ is an irreducible ${\mbox{sl}}(2,\C)$-module and formula 
(\ref{sl2c1}) applies. We have
$$
DU(x_{u}')=\alpha_u^2 x_{u}'=(u+1-k)(n-k-u)x_{u}'$$
and thus $\alpha_u = \sqrt{(u+1-k)(n-k-u)}$.

Using a little bit more representation theory 
we can give a characterization of $O(n)$ among all SJB's satisfying Theorem
\ref{mt2}.
Consider an irreducible $S_n$-module $V$. By the branching rule 
the decomposition of $V$ into irreducible $S_{n-1}$-modules is multiplicity
free and is therefore
canonical. Each of these modules, in turn, decompose canonically into
irreducible $S_{n-2}$-modules. Iterating this construction we get a
canonical decomposition of $V$ into irreducible $S_1$-modules, i.e., one
dimensional subspaces. Thus, there is a canonical basis of $V$, determined
upto scalars, and  called the
{\em Gelfand-Tsetlin or Young basis (GZ-basis)} (see {\bf \cite{vo}}). Note that the GZ-basis
is orthogonal wrt the (unique upto scalars) $S_n$-invariant inner product on
$V$. We now observe  the following:

(i) If $f:V\rar W$ is a $S_n$-linear isomorphism between
irreducibles $V,W$ then the GZ-basis of $V$ goes to the GZ-basis of $W$. 

(ii) Let $V$ be a $S_n$-module whose decomposition into irreducibles is
multiplicity free. By the GZ-basis of $V$ we mean the union of the GZ-bases of
the various irreducibles occuring in the (canonical) decomposition of $V$
into irreducibles. Then the GZ-basis of $V$ is orthogonal wrt any
$S_n$-invariant inner product on $V$.

Now
consider the $S_n$ action on $V(B(n))$. Since $U$ is $S_n$-linear, the
action is multiplicity free on $V(B(n)_k)$, for all $k$, and there exists a
SJB of $V(B(n))$, it follows from points (i) and (ii) above that there is a
canonically defined
(upto scalars) orthogonal SJB of $V(B(n))$ that consists of the 
union of the GZ-bases of $V(B(n)_k)$, $0\leq k\leq n$. We call this basis
the {\em symmetric Gelfand-Tsetlin basis} of $V(B(n))$.  
We prove the following result in Section 4.

\bt \label{mt3} 
The SJB $O(n)$, produced by the linear BTK algorithm when applied to the
Boolean algebra $B(n)$, is the symmetric Gelfand-Tsetlin basis of $V(B(n))$.
\et
In Example \ref{be} in Section 3 we write down the symmetric Gelfand-Tsetlin
bases of $V(B(n))$ for $n=1,2,3,4$. 

{\textcolor{black}{\bf \section { Terwilliger algebra of the binary Hamming
scheme }}}

The {\em Terwilliger algebra of the binary Hamming scheme}, denoted $\T_n$,
is defined to be the commutant of the $S_n$ action on $B(n)$, i.e., $\T_n =
\mbox{End}_{S_n}(V(B(n)))$ (in this definition the order structure on $B(n)$
is irrelevant). In this section we explicitly block diagonalize $\T_n$. 
It is convenient to think of
$B(n)$ as the poset of subsets of the set $\{1,2,\ldots ,n\}$.

Being the commutant of a finite group action
 $\T_n$ is a  $C^*$-algebra. Note that $\T_n$ is noncommutative (since, for
example, the trivial representation occurs at every rank and thus more than
once).
Let us first describe $\T_n$ in matrix terms. We represent
elements of $\mbox{End}(V(B(n)))$ (in the standard basis) as $B(n)\times
B(n)$ matrices (we think of elements of $V(B(n))$ as column vectors with
coordinates indexed by $B(n)$). For
$X,Y\in B(n)$, the entry in row $X$, column $Y$ of a matrix $M$ will be denoted
$M(X,Y)$. The matrix corresponding to $f\in \mbox{End}(V(B(n)))$ is denoted
$M_f$.
\bl \label{l1}
Let $f:V(B(n))\rar V(B(n))$ be a linear map. Then
$f$  is $S_n$-linear if and only if
$$M_f(X,Y)=M_f(\pi(X),\pi(Y)),\mbox{ for all } X,Y\in B(n),\;\pi\in S_n.$$
\el
\pf 
(only if) For $Y\in B(n)$  we have 
$f(Y) = \sum_{X\in B(n)} M_f(X,Y) X$.
Since $f$ is $S_n$-linear we now have, for $\pi \in S_n$,
$$\sum_{X\in B(n)} M_f(X,\pi(Y)) X = f(\pi(Y)) = \pi(f(Y))= 
\sum_{X\in B(n)} M_f(X,Y) \pi(X).$$
It follows that 
$M_f(X,\pi(Y))=M_f(\pi^{-1}(X),Y)$. Thus we have
$M_f(\pi(X),\pi(Y))=M_f(X,Y)$.

\noi (if) Similar to the only if part. $\Box$

Define $\A_n$ to be the set of all $B(n)\times B(n)$ complex matrices $M$
satisfying $M(X,Y)=M(\pi(X),\pi(Y))$, for all $X,Y\in B(n),\;\pi\in S_n$.
It follows from Lemma \ref{l1} that $\A_n$ is a $C^*$-algebra of matrices 
isomorphic to $\T_n$. We can easily determine its dimension. For nonnegative
integers $i,j,t$ let $M^t_{i,j}$ be the $B(n)\times B(n)$ matrix given by
$$M^t_{i,j}(X,Y) = \left\{ \ba{ll}
                            1 & \mbox{if }|X|=i,\;|Y|=j,\;|X\cap Y|=t \\
                            0 & \mbox{otherwise}
                            \ea
                    \right.
$$
Given $(X,Y),(X',Y')\in B(n)\times B(n)$, there exists $\pi\in S_n$ with
$\pi(X)=X',\;\pi(Y)=Y'$ if and only if $|X|=|X'|,\;|Y|=|Y'|$, and
$|X\cap Y|=|X'\cap Y'|$. It follows that 
$$\{M^t_{i,j}\;|\; i-t+t+j-t \leq n,\;i-t,t,j-t \geq 0 \}$$
is a basis of $\A_n$ and its cardinality is $\binom{n+3}{3}$.

It follows from general $C^*$-algebra theory that there exists a
{\em{block diagonalization}} of $\A_n$, i.e.,  there exists a 
$B(n)\times S$ unitary matrix $N(n)$, for some index set $S$ of cardinality
$2^n$,
and positive integers
$p_0,q_0,\ldots ,p_m,q_m$ such that ${N(n)}^*\A_nN(n)$ is equal to the 
set of all $S\times S$
block-diagonal matrices
\beq \label{bd1} &
\left( \ba{cccc}  C_0 & 0 & \ldots &0\\
            0   & C_1 &\ldots & 0\\
            \vdots & \vdots &\ddots &\vdots\\
            0 & 0 & \ldots & C_m \ea \right)
&
\eeq
where each $C_k$ is a block-diagonal matrix with $q_k$ repeated, identical
blocks of order $p_k$
\beq \label{bd2} &
C_k = \left( \ba{cccc}  B_k & 0 & \ldots &0\\
            0   & B_k &\ldots & 0\\
            \vdots & \vdots &\ddots &\vdots\\
            0 & 0 & \ldots & B_k \ea \right)
&
\eeq
Thus $p_0^2 + \cdots + p_m^2 = \mbox{dim}(\A_n)$ and $p_0q_0 + \cdots +
p_mq_m = 2^n$. The numbers $p_0,q_0,\ldots ,p_m,q_m$ and $m$ are uniquely
determined (upto permutation of the indices) by $\A_n$.

By dropping duplicate blocks we get a positive semidefiniteness 
 preserving $C^*$-algebra isomorphism (below $\mbox{Mat}(n\times n)$ denotes
the algebra of complex $n\times n$ matrices)
$$\Phi: \A_n \cong \bigoplus_{k=0}^m \mbox{Mat}(p_k\times p_k).$$

In an 
{\em explicit block diagonalization} we need to know this isomorphism
explicitly, i.e., we need to know  the entries in the
image of $M_{i,j}^t$. In {\bf \cite{s}}, an explicit block
diagonalization was determined. We now show that this result follows from
Theorem \ref{mt2}. 
The present proof also yields an explicit, canonical
(real) unitary matrix $N(n)$ achieving the isomorphism $\Phi$. 

The first step is a binomial inversion argument. Fix $i,j\in\{0,\ldots
,n\}$. 
Then we have
$$M_{i,t}^tM_{t,j}^t = \sum_{u=0}^n {u\choose t}M_{i,j}^u,\;\;\;\;\;
t=0,\ldots ,n,$$
since the entry of the lhs in row $X$, col $Y$ with $|X|=i, |Y|=j$ is equal
to the number of common subsets of $X$ and $Y$ of size $t$. 
Apply binomial inversion to get
\beq \label{bi}
M_{i,j}^t= \sum_{u=0}^n (-1)^{u-t}{u\choose
t}M_{i,u}^uM_{u,j}^u,\;\;\;\;\;t=0,\ldots ,n.
\eeq
Since $M_{u,j}^u = (M_{j,u}^u)^t$ and it will turn out that $N(n)$ can be
taken to be real (see the definition below) it follows that
\beq 
\Phi(M_{i,j}^t)= \sum_{u=0}^n (-1)^{u-t}{u\choose
t}\Phi(M_{i,u}^u)\Phi(M_{j,u}^u)^t,\;\;\;\;\;t=0,\ldots ,n,
\eeq
and hence all the images under $\Phi$ can be calculated by 
knowing the images $\Phi(M_{i,u}^u)$. 


For the second step we use Theorem \ref{mt2} whose notation we
preserve. For the rest of this section
set $m=\lfloor n/2 \rfloor$, and $p_k
=n-2k+1,\; q_k = {n\choose k}-{n\choose {k-1}},\; k=0,\ldots ,m$. 
Note that
\beq \label{di}
\sum_{k=0}^m p_k^2 = {{n+3}\choose 3},
\eeq
since both sides are polynomials in $l$ (treating the cases $n=2l$ and
$n=2l+1$ separately) of degree 3 and agree for $l=0,1,2,3$.

For
$0\leq k \leq m$, $O(n)$ will contain $q_k$ symmetric Jordan chains, each
containing $p_k$ vectors, starting 
at rank $k$ and ending at rank $n-k$. We can formalize this as follows: define the
finite set $S=\{(k,b,i)\;|\;0\leq k \leq m ,\; 1\leq b \leq q_k ,\; k \leq i \leq
n-k\}$. For each $0\leq k \leq m$, fix some linear ordering of the $q_k$
Jordan chains of $O(n)$ going from rank $k$ to rank $n-k$. 
Then there is a bijection ${\cal B}: O(n) \rar S$ defined as
follows: let $v\in O(n)$. Then ${\cal B}(v) = (k,b,i)$, where $i=r(v)$
and $v$ occurs on the $b$th
symmetric Jordan chain going from rank $k$ to rank $n-k$ (there are unique such
$k,b$). 
Linearly order
$S$ as follows: $(k,b,i) <_{\ell} (k',b',i')$ iff $k<k'$ or $k=k', b< b'$ or
$k=k',b=b', i< i'$. Form a $B(n)\times S$ matrix $N(n)$ as
follows: the columns of $N(n)$ are the normalized images $\frac{{\cal
B}^{-1}(s)}{\pr {\cal B}^{-1}(s) \pr}, s\in S$ listed in increasing order (of $<_\ell$).
By Theorem \ref{mt2}(i), $N(n)$ is unitary. 
Since the action of $M_{i,u}^u$ on
$V(B(n)_u)$ is $\frac{1}{(i-u)!}$ times the action of $U^{i-u}$ on
$V(B(n)_u)$, it follows  
by Theorem \ref{mt2}(ii) and identities (\ref{bi}), (\ref{di}) above 
that conjugating by
$N(n)$ provides a block diagonalization of $\A_n$ of the
form (\ref{bd1}), (\ref{bd2}) above. Set $\Phi$ equal to conjugation by
$N(n)$ followed by dropping duplicate blocks. To
calculate the images under $\Phi$ we shall now use part (iii) of Theorem
\ref{mt2}.

For $i,j,k,t \in \{0,\ldots ,n\}$ define
\beqn
\beta_{i,j,k}^t & = &\sum_{u=0}^n (-1)^{u-t} {u\choose t}{{n-2k}\choose
{u-k}}{{n-k-u}\choose {i-u}}{{n-k-u}\choose {j-u}}.
\eeqn
For $0\leq k \leq m$ and $k\leq i,j \leq n-k$, define $E_{i,j,k}$ to be the
$p_k \times p_k$ matrix, with rows and columns indexed by $\{k,k+1,\ldots
,n-k\}$, and with entry in row $i$ and column $j$ equal to 1 and all other
entries 0.
\bt  \label{scht}
{\em (}Schrijver {\bf\cite{s}}{\em )}
Let $i,j,t\in \{0,\ldots ,n\}$. Write
$$\Phi (M_{i,j}^t) = (N_0,\ldots ,N_m),$$
where, for
$k=0,\ldots , m$, the rows and columns of $N_k$
are indexed by $\{k,k+1,\ldots ,n-k\}$. 
Then, for $0\leq k \leq m$,
$$N_k = \left\{ \ba{ll}
                {{n-2k}\choose{i-k}}^{-\frac{1}{2}}
                {{n-2k}\choose{j-k}}^{-\frac{1}{2}}
                \beta_{i,j,k}^t
                E_{i,j,k} & \mbox{if } k\leq i,j \leq n-k \\
                  0    & \mbox{otherwise}
                 \ea
         \right.
$$ 
\et
\pf Fix $0\leq k \leq m$. If both $i,j$ are not elements of $\{k,\ldots ,
n-k\}$ then clearly $N_k = 0$. So we may assume $k\leq i,j \leq n-k$.
Clearly, $N_k = \lambda E_{i,j,k}$ for some $\lambda$. We now find
$\lambda = N_k(i,j)$.

Let $u\in \{0,\ldots ,n\}$.
Write 
$\Phi (M_{i,u}^u) = (A_0^u,\ldots ,A_m^u)$.
We claim that 
$$A_k^u = \left\{ \ba{ll}
                {{n-k-u}\choose{i-u}}
                {{n-2k}\choose{u-k}}^{\frac{1}{2}}
                {{n-2k}\choose{i-k}}^{-\frac{1}{2}}
                E_{i,u,k} & \mbox{if } k\leq u \leq n-k \\
                  0    & \mbox{otherwise}
                 \ea
         \right.
$$ 
The otherwise part of the claim is clear. If $k\leq u \leq n-k$ and $i < u$
then we have $A_k^u = 0$. This also follows from the rhs since the binomial
coefficient $a\choose b$ is 0 for $b< 0$. So we may assume that $k\leq u
\leq n-k$ and $i\geq u$. Clearly, in this case we have $A_k^u = \alpha
E_{i,u,k}$, for some $\alpha$. We now determine $\alpha= A^u_k(i,u)$. We
have using Theorem \ref{mt2}(iii) 
$$
A_k^u(i,u) =  \frac{\prod_{w=u}^{i-1}\left\{
(n-k-w){{n-2k}\choose{w-k}}^{\frac{1}{2}}
                {{n-2k}\choose{w+1-k}}^{-\frac{1}{2}}\right\}}{(i-u)!}
=
{{n-k-u}\choose{i-u}}{{n-2k}\choose{u-k}}^{\frac{1}{2}}
                {{n-2k}\choose{i-k}}^{-\frac{1}{2}}
$$

Similarly, if we write 
$\Phi (M_{u,j}^u) = (B_0^u,\ldots ,B_m^u)$
then, since $M_{u,j}^u = (M_{j,u}^u)^t$, we have
$$B_k^u = \left\{ \ba{ll}
                {{n-k-u}\choose{j-u}}
                {{n-2k}\choose{u-k}}^{\frac{1}{2}}
                {{n-2k}\choose{j-k}}^{-\frac{1}{2}}
                E_{u,j,k} & \mbox{if } k\leq u \leq n-k \\
                  0    & \mbox{otherwise}
                 \ea
         \right.
$$
It now follows from (\ref{bi}) that 
$
N_k  =  \sum_{u=0}^n (-1)^{u-t} {u\choose t} A_k^u B_k^u  
 =  \sum_{u=k}^{n-k} (-1)^{u-t} {u\choose t} A_k^u B_k^u. 
$
Thus  
\beqn
\lefteqn{N_k(i,j)}\\
 &=& \sum_{u=k}^{n-k} (-1)^{u-t} {u\choose t}  
\left\{ \sum_{l=k}^{n-k} A_k^u(i,l) B_k^u(l,j)\right\}\\  
 &=&  \sum_{u=k}^{n-k} (-1)^{u-t} {u\choose t}  A_k^u(i,u) B_k^u(u,j)\\
 &=& \sum_{u=k}^{n-k} (-1)^{u-t} {u\choose t}
{{n-k-u}\choose{i-u}}{{n-2k}\choose{u-k}}^{\frac{1}{2}}
                {{n-2k}\choose{i-k}}^{-\frac{1}{2}}
{{n-k-u}\choose{j-u}}{{n-2k}\choose{u-k}}^{\frac{1}{2}}
                {{n-2k}\choose{j-k}}^{-\frac{1}{2}}\\
& = & {{n-2k}\choose{i-k}}^{-\frac{1}{2}}
      {{n-2k}\choose{j-k}}^{-\frac{1}{2}}
      \left\{\sum_{u=0}^{n} (-1)^{u-t} {u\choose t}
{{n-k-u}\choose{i-u}}
{{n-k-u}\choose{j-u}}{{n-2k}\choose{u-k}}\right\}.
\eeqn

{\textcolor{black} {\bf \section {The linear BTK algorithm }}}

In this section we present the linear analog of the BTK algorithm for
constructing a SJB of $V(M(n,k_1,\ldots ,k_n))$. Though we do not recall
here the BTK algorithm for constructing a SCD of $M(n,k_1,\ldots ,k_n)$ 
(see {\bf \cite{btk,a,e}}),
readers familiar with that method will easily recognize the present
algorithm as its linear analog.

The basic building block of the linear BTK algorithm is an inductive method
for constructing a SJB of $V(M(2,p,q))$. If $p \mbox{ or } q = 0$,
then $M(2,p,q)$ is order isomorphic to a chain and the characteristic
vectors of the elements of the chain form a SJB .
For positive $p,q$ we shall now reduce the problem of constructing a SJB of
$V(M(2,p,q))$ to that of constructing a SJB of $V(M(2,p - 1,q - 1))$. We begin
with the following elementary lemma on determinants. 

\bl \label{dl}
Let $N=(a_{i,j})$ be a $n\times n$ real matrix, $n\geq 2$. Suppose that 
 
\noi (i) $a_{i,1} > 0$, for $i\in \{1,\ldots ,n\}$, i.e., 
the first column contains positive entries. 

\noi (ii) For $j\in \{2,\ldots ,n\}$, $a_{j,j} > 0$, $a_{j-1,j} < 0$, and
all other entries in column $j$ are 0.

\noi Then $\mbox{det}(N) > 0$.
\el
\pf By induction on $n$. The assertion is clear for $n=2$. Now assume that
$n>2$. Let $N[i,j]$ denote the $(n-1)\times (n-1)$ matrix obtained from $N$
by deleting row $i$ and column $j$. Expanding $\mbox{det}(N)$ by the first
row we get
$$\mbox{det}(N) = a_{1,1}\mbox{det}(N[1,1]) - a_{1,2}\mbox{det}(N[1,2])$$
Now $a_{1,1}>0$, $a_{1,2}<0$, $N[1,1]$ is upper triangular with positive
diagonal entries, and $\mbox{det}(N[1,2]) > 0$ (by induction hypothesis).
The result follows.$\Box$

Let $p,q$ be positive and set $P=M(2,p,q)$, $W=V(P)$ with up
operator $U$. Let $r$ denote the rank function of $P$. We have $\mbox{dim
}W=(p+1)(q+1)$. The action of $U$ on the standard basis of $W$ is given as
follows: 
for $0\leq i \leq p$, $0\leq j \leq q$
\beqn
U((i,j)) &=& \left\{ \ba{ll}
      (i+1,j) + (i,j+1) & \mbox{if $i<p,j<q$}\\ 
      (i+1,j)  & \mbox{if $i<p,j=q$}\\ 
      (i,j+1) & \mbox{if $i=p,j<q$}\\ 
      0 & \mbox{if $i=p,j=q$} 
                   \ea \right.
\eeqn
Consider the following symmetric Jordan chain in $W$ generated by
$v(0)=(0,0)$:
$$(v(0),v(1),v(2),\ldots,v(p+q)),$$
where, for $0\leq k \leq p+q$,
\beq \label{gc1} v(k)=U^k((0,0))=\sum_{i,j}{k\choose i}(i,j),\eeq
the sum being over all $0\leq i \leq p$, $0\leq j \leq q$ with $i+j=k$.

The following result is basic to our inductive approach.
\bt \label{it}
Define homogeneous vectors in $W$ as follows:
\beq\label{gc2}
v(i,j)= (p-i)(i,j) - (q-j+1)(i+1,j-1),\;\;0\leq i \leq p-1,\;1\leq j \leq
q.
\eeq
Then

\noi (i) $v(i,j)$ is nonzero and $r(v(i,j)) = i+j,\;\;0\leq i \leq
p-1,\;1\leq j \leq q.$ 

\noi (ii) $\{v(k) \;|\; 0\leq k \leq p+q \} \cup \{v(i,j) \;|\; \;\;0\leq i \leq
p-1,\;1\leq j \leq q \}$ is a basis of $W$.

\noi (iii) For $0\leq i \leq p-1$, $1\leq j \leq q$ we have
\beqn
U(v(i,j)) &=& \left\{ \ba{ll}
      v(i+1,j) + v(i,j+1) & \mbox{if $i<p-1,j<q$}\\ 
      v(i+1,j)  & \mbox{if $i<p-1,j=q$}\\ 
      v(i,j+1) & \mbox{if $i=p-1,j<q$}\\ 
      0 & \mbox{if $i=p-1,j=q$} 
                   \ea \right.
\eeqn
Thus, the action of $U$ on the $v(i,j)$ is isomorphic to the action
of the up operator on the standard basis of $V(M(2,p-1,q-1))$, except
that the map $(i,j) \mapsto v(i,j+1),\;\;(i,j)\in M(2,p-1,q-1)$ shifts ranks
by one $($since $r(v(i,j+1))=i+j+1$$)$.
\et
\pf  (i) This is clear.

\noi (ii) For $0\leq k \leq p+q$ define
$$X_k = \{ v(k)\} \cup \{v(i,j)\;:\;0\leq i\leq p-1,\; 1\leq j \leq q,\; i+j =
k\}$$
The map $\phi_k : X_k \rar M(2,p,q)_k$ given by
$\phi_k(v(i,j))=(i,j)$ and 
$$\phi_k (v(k)) =\left\{ \ba{ll}
                         (k,0) & \mbox{if } k<p \\
                         (p,k-p) & \mbox{if } k\geq p
                         \ea \right. $$
is clearly a bijection. It is enough to show that $X_k$ is a basis of
$V(M(2,p,q)_k)$. Linearly order the elements of $M(2,p,q)_k$ using reverse
lexicographic order $<_r$: $(i,j) <_r (i',j')$ if and only if $i>i'$.
Transfer this order to $X_k$ via $\phi_k^{-1}$. Consider the $M(2,p,q)_k
\times X_k$ matrix $N$, with rows and columns listed in the order $<_r$, and
whose columns are the coordinate vectors of elements of $X_k$ in the standard basis
$M(2,p,q)_k$  of $V(M(2,p,q)_k)$. Equation (\ref{gc1}) shows that hypothesis
(i) of Lemma \ref{dl} is satisfied and equation (\ref{gc2}) shows that
hypothesis (ii) of Lemma \ref{dl} is satisfied. The result now follows from
Lemma \ref{dl}.   

\noi (iii) We check the first case. The other cases are similar. 
Let $i<p-1, j<q$. Then
\beqn \lefteqn{U(v(i,j))}\\ 
&=& U((p-i)(i,j) - (q-j+1)(i+1,j-1)) \\
                &=& (p-i)((i+1,j)+(i,j+1)) - (q-j+1)((i+2,j-1) + (i+1,j)) \\
                &=& (p-i-1)(i+1,j) - (q-j+1)(i+2,j-1) + (p-i)(i,j+1) -
(q-j)(i+1,j)\\
                &=& v(i+1,j) + v(i,j+1),
\eeqn 
completing the proof.$\Box$

Theorem \ref{it}, whose notation we preserve, gives the following inductive
method for constructing a SJB of $V(M(2,p,q))$: if $p$ or $q$ equals 0,
the chain $M(2,p,q)$ itself gives a SJB. 
Now suppose $p,q > 1$. Set $v(0)=(0,0)$ and form the
symmetric Jordan chain $C=(v(0),v(1),\ldots ,v(p+q))$, where $v(k)$ is given
by (\ref{gc1}). Take the (inductively constructed) 
SJB of $V(M(2,p-1,q-1))$ and (using parts (ii) and  (iii) of
Theorem \ref{it}) transfer each Jordan
chain in this SJB to a Jordan chain in $V(M(2,p,q))$ via the map $(i,j)\mapsto
v(i,j+1)$. Since $r(v(0,1))=1$ and $r(v(p-1,q))=p+q-1$, 
each such transferred Jordan chain is symmetric in
$V(M(2,p,q))$ and part (ii) of Theorem \ref{it} now shows that 
the collection of these chains together with $C$ gives a SJB
of $V(M(2,p,q))$. 

\bex \label{ex}
{\em 

\noi (i) Here we work out the SJB of $V(M(2,2,2))$ produced by the algorithm
above. The symmetric Jordan chain generated by $(0,0)$ is given by
\beq \label{ch1}
((0,0)\,,\,(1,0)+(0,1)\,,\,(2,0)+2(1,1)+(0,2)\,,\,3(2,1)+3(1,2)\,,\,6(2,2)).
\eeq
The $v(i,j)$ are given by
\beqn
v(0,1)=2(0,1)-2(1,0)&&v(1,1)=(1,1)-2(2,0) \\
v(0,2)=2(0,2)-(1,1)&&v(1,2)=(1,2)-(2,1)
\eeqn 
The SJB of $V(M(2,1,1))$ is given by the following two chains
\beqn
&((0,0)\,,\,(1,0)+(0,1)\,,\,2(1,1))&\\
&((0,1)-(1,0))&
\eeqn
Transferring these chains to $V(M(2,2,2))$ via the map $(i,j)\mapsto
v(i,j+1)$ gives the following two chains
\beq \label{ch2}
&(v(0,1)\,,\,v(1,1)+v(0,2)\,,\,2v(1,2))&\\
\label{ch3}
&(v(0,2)-v(1,1))&
\eeq
Chains (\ref{ch1}, \ref{ch2}, \ref{ch3}) give a SJB of $V(M(2,2,2))$.

(ii) The procedure above is especially simple when, say
$q = 1$, as in this case the recursion stops right at the first stage.
For later reference we spell out this case in detail. Consider $M(2,n,1)$.
Define a 1-1 linear map $V(M(2,n,0)) \rar V(M(2,n,1))$ by $(i,0)\mapsto
(i,1),\;(i,0)\in M(2,n,0)$. For $v\in V(M(2,n,0))$, we denote the image of
$v$ under this map by $\ol{v}$.

Let $(x_0,x_1,\ldots ,x_n)$, where $x_i = (i,0)$ be the SJB of $V(M(2,n,0))$.
Set $x_{-1}=x_{n+1}=0$. We now consider two cases:

(a) $n=0$ : In this case 
\beq \label{c0}
&(x_0,\ol{x_0}).&
\eeq
is the SJB of $V(M(2,n,1))$ produced by Theorem \ref{it}.

(b) $n\geq 1$ : The symmetric Jordan chain in $V(M(2,n,1))$ generated by
$(0,0)$ can be written as (using (\ref{gc1}))
\beq \label{c1}
(y_0,y_1,\ldots ,y_{n+1}),&&
\mbox{ where }\;
y_l = x_l + l\, \ol{x_{l-1}},\;\;0\leq l \leq n+1.\eeq
The SJB of $V(M(2,n,1))$ produced by Theorem \ref{it} is given by (\ref{c1})
and the following symmetric Jordan chain:
\beq \label{c3}
(z_1,\ldots ,z_n),&&
\mbox{ where }\;
z_l = (n-l+1)\,\ol{x_{l-1}} -x_l,\;\;1\leq l \leq n.\eeq

}\eex

We can now give the full linear BTK algorithm which reduces the general case to the
$n=2$ case.

\pf {\bf of Theorem \ref{mt1}} The proof is by induction on $n$, the case $n=1$
being clear and the case $n=2$ established above. Let $P=M(n,k_1,\ldots
,k_n)$, $n\geq 3$ and set $V=V(P)$. Denote the rank function of $P$ by $r$. 
Define induced subposets $P(j)$ of $P$ by
$$ P(j) = \{(a_1,\ldots ,a_n)\in P\;:\; a_n = j \},\;\;\;0\leq j 
\leq k_n,$$
and set $V(j) = V(P(j))$. Let $U$ denote the up operator on $V$ and
$U_j$ denote the up operator on $V(j)$, $0 \leq j \leq k_n$. Note
that all the $P(j)$ are order isomorphic. We have
\beq \label{p6}
V &=& V(0)\oplus V(1) \oplus \cdots \oplus V(k_n).
\eeq
For $0\leq j < k_n$ define linear isomorphisms $R_j:V(j) \rar V(j+1)$ by
$$R_j((a_1,\ldots ,a_{n-1},j))=(a_1,\ldots ,a_{n-1},j+1),\;\;(a_1,\ldots
,a_{n-1},j)\in P(j)$$
Put $V(k_n+1)=\{0\}$ and define $R_{k_n}$ to be the zero map. 

For $v\in V(j),\; 0\leq j \leq k_n$ we have
\beq
\label{p1}
U(v)&=&U_j(v) + R_j(v)\\
\label{p2}
U_{j+1}R_j(v)&=&R_jU_j(v)
\eeq
By induction there is a SJB of $V(0)$ (wrt $U_0$). Let $t$ denote the number
of symmetric Jordan chains in this SJB, with the $m^{th}$ chain denoted by
$$S(0,m)=(v(0,0,m),v(1,0,m),\ldots ,v(l_m,0,m)),\;\;\;1\leq m \leq t,$$
where $l_m\geq 0$ is the length of the $m^{th}$ symmetric Jordan chain. Thus  
$$\ba{ll}
      r(v(0,0,m)) + r(v(l_m,0,m)) = k_1+\cdots +k_{n-1} & \mbox{if
$l_m>0$}\\ 
      2r(v(0,0,m)) = k_1+\cdots +k_{n-1} & \mbox{if $l_m=0$} 
                   \ea 
$$
Define subspaces $X(0,m)\subseteq V(0)$ by
$$ X(0,m) = \mbox{Span }S(0,m)
             = \mbox{Span }\{v(0,0,m),\ldots ,v(l_m,0,m)\},\;\;1\leq m
\leq t
$$
We have ${\ds{V(0) = X(0,1)\oplus \cdots \oplus X(0,t)}}$. Note that
$\mbox{dim }X(0,m)=l_m+1$.

\noi For $1\leq m \leq t,\; 1\leq j \leq k_n,\; 0\leq i \leq l_m$ define, by
induction on $j$ (starting with $j=1$),
\beq
\label{p3}
v(i,j,m)&=&R_{j-1}(v(i,j-1,m))\\
\nonumber
S(j,m)&=& (v(0,j,m),v(1,j,m),\ldots ,v(l_m,j,m)) \\ \label{p10}
X(j,m) &=& \mbox{Span }S(j,m) = \mbox{Span }\{v(0,j,m),\ldots ,v(l_m,j,m)\}
\eeq
%
Since the $R_j$, $0\leq j < k_n$ are isomorphisms we have
\beq \label{p8}
&& \{v(0,j,m),\ldots ,v(l_m,j,m)\} \mbox{ is independent
},\;0\leq j \leq k_n,1\leq m \leq t. \\
\label{p7} &&
V(j) = X(j,1)\oplus X(j,2)\oplus \cdots \oplus X(j,t),\;\;0\leq j \leq
k_n.
\eeq
For $1\leq m \leq t,\;1\leq j \leq k_n$ we see from (\ref{p8}) and 
the following inductive calculation on $j$ (using (\ref{p2})) that
each $S(j,m)$ is a graded Jordan chain in $V(j)$ (below we take
$v(l_m+1,0,m)=0$, for all $m$) 
\beq \nonumber
U_j(v(i,j,m))&=&U_jR_{j-1}(v(i,j-1,m))\\ \nonumber
             &=&R_{j-1}U_{j-1}(v(i,j-1,m))\\ \nonumber
             &=&R_{j-1}(v(i+1,j-1,m))\\ \label{p9}
             &=&v(i+1,j,m)
\eeq
For $m=1,\ldots ,t$
define subspaces $Y(m)\subseteq V$ by
\beq \label{p11} Y(m) &=&
 X(0,m)\oplus X(1,m) \oplus \cdots \oplus X(k_n,m). \eeq
We have from (\ref{p6}) and (\ref{p7}) that
\beq \label{p12}
V &=& Y(1)\oplus \cdots \oplus Y(t).
\eeq
It follows from (\ref{p10}), (\ref{p8}), and (\ref{p11}) that the set
$$B(m)=\{v(i,j,m)\;:\;0\leq i \leq l_m,\;0\leq j\leq k_n \},\;\;1\leq m
\leq t$$
is a basis of $Y(m)$. Note that
\beq \label{p13}
r(v(i,j,m)) &=& r(v(0,0,m)) + i + j,\;\;1\leq m \leq t,\; 
0\leq i \leq l_m,\;0\leq j\leq k_n.\eeq

Fix $1\leq m \leq t$. Using (\ref{p1}), (\ref{p3}), and (\ref{p9}) we see 
that the action of $U$ on the basis $B(m)$ is given by
\beqn
U(v(i,j,m)) &=& \left\{ \ba{ll}
      v(i+1,j,m) + v(i,j+1,m) & \mbox{if $i<l_m,j<k_n$}\\ 
      v(i+1,j,m)  & \mbox{if $i<l_m,j=k_n$}\\ 
      v(i,j+1,m) & \mbox{if $i=l_m,j<k_n$}\\ 
      0 & \mbox{if $i=l_m,j=k_n$} 
                   \ea \right.
\eeqn
So this action is isomorphic to the action of the up operator on the
standard basis of $V(M(2,l_m,k_n))$, except for the shift of rank given by
(\ref{p13}). We can now use the  algorithm of Theorem \ref{it} to construct
an SJB of $Y(m)$ wrt $U$.
Since 
$$r(v(0,0,m))+r(v(l_m,k_n,m)) = r(v(0,0,m))+r(v(l_m,0,m)) + k_n = k_1 +
\cdots + k_n,$$
each graded Jordan chain in this SJB is symmetric in $V$. From (\ref{p12})
it follows that the union of the SJB's of $Y(m)$ gives an SJB of $V$. That completes the
proof.$\Box$

\pf {\bf of Theorem \ref{mt2}} The proof is by induction of $n$. The result is
clear for $n=1$. We write $B(n)$ as $M(n,k_1,\ldots ,k_n)$, where
$k_1=\cdots =k_n=1$.

Consider $V=V(M(n+1,k_1,\ldots ,k_{n+1}))$, with $k_i = 1$ for all $i$. 
We preserve the notation of the proof of Theorem \ref{mt1}. Let there be $t$
symmetric Jordan chains in the SJB of $V(0) = V(B(n))$. For $v\in V(0)$, we
denote $R_0(v) = \ol{v}$ (agreeing with the notation in Example \ref{ex}(ii)). 

We have 
\beqn
V&=& V(0) \oplus V(1)\\
V(0)&=&X(0,1)\oplus X(0,2) \oplus \cdots \oplus X(0,t) \\
V(1)&=&X(1,1)\oplus X(1,2) \oplus \cdots \oplus X(1,t) 
\eeqn
Since the inner product is standard we have
\beq \label{ipi}
\langle u,v \rangle = \langle \ol{u} ,\ol{v} \rangle ,\; \langle u, \ol{v}
\rangle  = 0,\;\;u,v\in V(0).\eeq
It follows that $V(0)$ is orthogonal to $V(1)$. The subspaces 
$X(0,1),\ldots ,X(0,t)$ are mutually orthogonal, by the induction
hypothesis, and thus  
\beqn V&=&Y(1)\oplus \cdots \oplus Y(t) \eeqn
is an orthogonal decomposition of $V$, where
\beqn Y(m)&=&X(0,m)\oplus X(1,m),\;\;1\leq m \leq t. \eeqn
Fix $1\leq m \leq t$. The subspaces $X(0,m)$ and $X(1,m)$ are orthogonal but
Theorem \ref{it} will produce  new symmetric Jordan chains from linear
combinations of vectors in $X(0,m)$ and $X(1,m)$. We now show that these too
are orthogonal and satisfy (\ref{mks}). This will prove the theorem.

Write the $m^{\mbox{th}}$ symmetric chain in $V(0)$ as 
\beq \label{ba1}
&(x_k,\ldots ,x_{n-k}),\;\;0\leq k \leq \lfloor n/2 \rfloor,&
\eeq
where $r(x_k) = k$. By induction hypothesis
\beq \label{i1}
\frac{\langle x_{u+1},x_{u+1}\rangle}{\langle x_u,x_u\rangle} &=& (u+1-k)(n-k-u),\;\;k\leq u < n-k.
\eeq
We now consider two cases.

(a) $k=n-k$ : It follows from (\ref{c0}) (after changing $n$ to $n-2k$ and 
shifting the rank by $k$) that the SJB 
of $Y(m)$ will consist of the single symmetric Jordan
chain 
\beq \label{ba2}
&(x_k, \ol{x_k}).&
\eeq 
We have
$$\frac{\langle\ol{x_k},\ol{x_k}\rangle}{\langle x_k,x_k\rangle} 
=\frac{\langle x_k,x_k\rangle}{\langle x_k,x_k\rangle}=1=(k+1-k)(n+1-k-k).$$

(b) $k< n-k$ : By (\ref{c1}, \ref{c3}) the SJB of $Y(m)$
will consist of the following two symmetric Jordan chains (after changing
$n$ to $n-2k$ and shifting the rank by $k$): 
\beq \label{ba3}
&(y_k,\ldots ,y_{n+1-k}), \mbox{ and }
(z_{k+1},\ldots ,z_{n-k}), \eeq
where
\beq \label{ba4} 
y_l &=& x_l + (l-k)\, \ol{x_{l-1}},\;\;k\leq l \leq n+1-k. \\\label{ba5}
z_l &=& (n-k-l+1)\,\ol{x_{l-1}} -x_l,\;\;k+1\leq l \leq n-k.\eeq
and $x_{k-1}=x_{n+1-k}=0$.

For $k+1\leq l \leq n-k$ we have from (\ref{ipi})
$$\langle y_l,z_l\rangle = (n-k-l+1)(l-k)\langle x_{l-1},x_{l-1}\rangle  - 
\langle x_l,x_l\rangle =0,$$
where the last step follows from (\ref{i1}) upon substituting $u=l-1$.
Thus 
$$\{y_k,\ldots ,y_{n+1-k},z_{k+1},\ldots ,z_{n-k}\}$$
is an orthogonal
basis of $Y(m)$.

We now check (\ref{mks}) for the case $n+1$. Let $k\leq u < n+1-k$. Note
that (\ref{i1}) holds for $u=n-k$ also. Also note that in
the following computation (in the
last but one step) we have used (\ref{i1}) when $u=k-1$ (in which case the rhs
is $0$ and the lhs is $\infty$). This is permissible here because of
the presence of the factor $(u-k)^2$. 

We have,
using (\ref{ipi}) and (\ref{i1}),
\beqn
\frac{\langle y_{u+1},y_{u+1}\rangle}{\langle y_u,y_u\rangle }&=&
\frac{\langle x_{u+1}+(u+1-k)\,\ol{x_u}\,,\,x_{u+1}+(u+1-k)\,\ol{x_u}\rangle}
{\langle x_{u}+(u-k)\,\ol{x_{u-1}}\,,\,x_{u}+(u-k)\,\ol{x_{u-1}}\rangle}\\
&=&\frac{\langle x_{u+1},x_{u+1}\rangle + (u+1-k)^2 \langle x_u,x_u\rangle}
{\langle x_{u},x_{u}\rangle + (u-k)^2 \langle x_{u-1},x_{u-1}\rangle}\\
&=&\frac{\frac{\langle x_{u+1},x_{u+1}\rangle}{\langle x_u,x_u\rangle} + (u+1-k)^2}
{1 + \frac{\langle x_{u-1},x_{u-1}\rangle}{\langle x_u,x_u\rangle}(u-k)^2}\\
&=& \frac{(u+1-k)(n-k-u) + (u+1-k)^2}{1+\frac{(u-k)^2}{(u-k)(n-k-u+1)}}\\
&=& (u+1-k)(n+1-k-u)
\eeqn
The calculation for $\frac{\langle z_{u+1},z_{u+1}\rangle}{\langle
z_{u},z_{u}\rangle}$ is similar and is
omitted.$\Box$

\bex \label{be}
{\em
In this example we work out the SJB's of $V(B(n))$, for $n=2,3,4$, starting
with the SJB of $V(B(1))$, using the formulas 
(\ref{ba1}, \ref{ba2}, \ref{ba3}, \ref{ba4}, \ref{ba5}) 
given in the proof of Theorem \ref{mt2}.
We write elements of $B(n)$ as subsets of $[n]=\{1,2,\ldots ,n\}$ rather than as
their characteristic vectors. Thus, for $X\subseteq [n]$, we have $R_0(X) =
\ol{X} = X \cup \{n+1\}$.

\noi (i) The SJB of $V(B(1))$ is given by
\beqn
&(\,\,\emptyset \,,\,\{1\}\,\,)&
\eeqn

\noi (ii) The SJB of $V(B(2))$ consists of 
\beqn
&(\,\,\emptyset \,,\, \{1\} + \{2\}\,,\, 2\{1,2\}\,\,)&\\
&(\,\,\{2\} - \{1\}\,\,)&
\eeqn

\noi (iii) The SJB of $V(B(3))$ consists of
\beqn
&(\,\,\emptyset \,,\, \{1\} + \{2\} + \{3\}\,,\,
2(\{1,2\}+\{1,3\}+\{2,3\})\,,\,6\{1,2,3\}\,\,)&\\
&(\,\,2\{3\}-\{1\}-\{2\}\,,\,\{1,3\}+\{2,3\}-2\{1,2\}\,\,)&\\
&(\,\,\{2\}-\{1\}\,,\,\{2,3\}-\{1,3\}\,\,)&
\eeqn

\noi (iv) The SJB of $V(B(4))$ consists of (some of the chains are split
across two lines)
\beqn
&&(\,\,\emptyset \,\,,\,\,\{1\}+\{2\}+\{3\}+\{4\}\,\,
,\,\, 2(\{1,2\}+\{1,3\}+\{1,4\}+\{2,3\}+\{2,4\}+\{3,4\})\,\,,\\
&&6(\{1,2,3\}+\{1,2,4\}+\{1,3,4\}+\{2,3,4\})\,\,,\,\,24\{1,2,3,4\}\,\,)\\
&&\\
&&(\,\,3\{4\} -
(\{1\}+\{2\}+\{3\})\,\,,
\,\,2(\{1,4\}+\{2,4\}+\{3,4\})-2(\{1,2\}+\{1,3\}+\{2,3\})\,\,,\\
&&\,\,2 (\{1,2,4\}+\{1,3,4\}+\{2,3,4\} )-6 \{1,2,3\} \,\,)\\
&&\\
&&(\,\,2\{3\}-(\{1\}+\{2\})\,\,,\,\,\{1,3\}
+\{2,3\}-2\{1,2\}+2\{3,4\}-(\{1,4\}+\{2,4\})\,\,,\\
&&\,\,
2(\{1,3,4\}+\{2,3,4\})-4\{1,2,4\}\,\,)\\
&&\\
&&(\,\,\{2\}-\{1\}\,\,,\,\,\{2,3\}-\{1,3\}+
\{2,4\}-\{1,4\}\,\,,\,\,2(\{2,3,4\}-\{1,3,4\})\,\,)\\
&&\\
&&(\,\,\{2,4\}-\{1,4\}-\{2,3\}+\{1,3\}\,\,)\\
&&\\
&&(\,\,2(\{3,4\}+\{1,2\})-(\{1,4\}+\{2,4\}+\{1,3\}+\{2,3\})\,\,)
\eeqn

It may be verified that each of the SJB's are orthogonal and satisfy
(\ref{mks}). 
}
\eex

{\textcolor{black}{\bf \section {Symmetric Gelfand-Tsetlin basis }}}

In this section we prove Theorem \ref{mt3}. The proof is an 
application of the Vershik-Okounkov {\bf\cite{vo}} theory of (complex) 
irreducible representations of the symmetric group. We first recall briefly
(without proofs) those points of the theory which we need: 

\noi (A) A direct elementary  argument is given to show that branching from
$S_n$ to $S_{n-1}$ is simple, i.e., multiplicity free. Once this is done we
have the canonically defined (upto scalars) Gelfand-Tsetlin basis (or GZ-basis) of an
irreducible $S_n$-module, as in the introduction. As stated there, the
GZ-basis of an irreducible representation $V$ is orthogonal wrt the unique
(upto scalars) $S_n$-invariant inner product on $V$.

\noi (B) Denote by $S_n^{\wedge}$ the set of
equivalence classes of finite dimensional complex irreducible
representations of $S_n$. Denote by $L^{\lambda}$ the irreducible $S_n$-module 
corresponding to $\lambda\in S_n^{\wedge}$.

We have identified a canonical basis, namely the GZ-basis,
in each irreducible representation of $S_n$. A natural question at this
point is to identify those elements of $\C[S_n]$ that act diagonally in this
basis (in every irreducible representation). In other words, consider the
algebra isomorphism
\beq\label{iso}
\C[S_n]&\cong&\bigoplus_{\lambda\in S_n^{\wedge}} \mbox{End}(L^{\lambda}),
\eeq
given by
$$\pi \mapsto ( L^{\lambda} \buildrel {\pi}\over \rightarrow L^{\lambda}\;:\;
\lambda
\in S_n^{\wedge}),\;\;\pi\in S_n.$$
Let
$\mbox{D}(L^{\lambda})$ consist of all operators on $L^{\lambda}$
diagonal in the GZ-basis of $L^{\lambda}$.
The question above  can now be stated as: what is the image under the isomorphism
(\ref{iso})
of the subalgebra $\bigoplus_{\lambda\in S_n^{\wedge}}
\mbox{D}(L^{\lambda})$ of
$\bigoplus_{\lambda\in S_n^{\wedge}}
\mbox{End}(L^{\lambda})$.

Let $Z_n$ denote the center of the algebra $\C[S_n]$ and set
 $GZ_n$ equal to the subalgebra of $\C[S_n]$ generated by $Z_1\cup Z_2\cup
 \cdots \cup Z_n$ (where we have the natural inclusions $S_1\subseteq S_2
\subseteq \cdots$). It is easy to see that $GZ_n$ is a commutative subalgebra of $\C[S_n]$.
It is called the {\em
Gelfand-Tsetlin algebra (GZ-algebra)} of the inductive family of group
algebras $\C[S_n]$. It is proved that
$GZ_n$ is the image of $\bigoplus_{\lambda\in S_n^{\wedge}}
\mbox{D}(L^{\lambda})$ under the isomorphism (\ref{iso}) above, i.e.,
$GZ_n$ consists of all elements of $\C[S_n]$ that act diagonally in
the GZ-basis in every irreducible representation of $S_n$. Thus
$GZ_n$ is a maximal commutative subalgebra of $\C[S_n]$ and its dimension is
equal to ${\sum_{\lambda\in S_n^{\wedge}} \mbox{ dim }L^{\lambda}}$.
By a GZ-vector $v$ we mean an element of the GZ-basis of $L^{\lambda}$, for
some $\lambda \in S_n^{\wedge}$. It follows that the GZ-vectors are the only
vectors that are eigenvectors for the action of every element of $GZ_n$. 
Moreover, any GZ-vector  
is uniquely determined by the eigenvalues
of the elements of $GZ_n$ on this vector.



\noi (C) For $i=1,2,\ldots ,n$ define $X_i = (1, i) + (2, i) + \cdots +
(i-1, i) \in \C[S_n]$. The $X_i$'s are called the
Young-Jucys-Murphy elements (YJM-elements) and it is shown that
 they generate $GZ_n$.
To a GZ-vector $v$  we associate the tuple $\alpha (v) =
(a_1,a_2,\ldots ,a_n)$,
where $a_i = \mbox{ eigenvalue of $X_i$  on } v$. We call $\alpha(v)$ the
{\em weight} of $v$ and we set
$$\spec(n) = \{\alpha (v)\;:\;v \mbox{ is a GZ-vector}\}.$$
It follows from step (B) above that, for GZ-vectors $u$ and $v$, $u=v$ iff
$\alpha(u) = \alpha(v)$ and thus 
$\# \spec(n)={\sum_{\lambda\in S_n^{\wedge}} \mbox{ dim }L^{\lambda}}$.
Given $\alpha \in \spec (n)$ we denote by $v_\alpha$ ($\in L^{\lambda}$ for
some unique $\lambda\in S_n^{\wedge}$) the GZ-vector with
weight $\alpha$.

There is a natural equivalence relation $\sim$ on $\spec (n)$: for $\alpha,
\beta \in \spec (n)$,
\beqn \alpha \sim \beta &\Leftrightarrow&
v_\alpha \mbox{ and } v_\beta \mbox{
belong to the same irreducible $S_n$-module $L^{\lambda}$ for some $\lambda
\in S_n^{\wedge}$}
\eeqn
Clearly we have
$\#(\spec (n)/\sim)=\# S_n^{\wedge}=\mbox{ number of partitions of }n$.

\noi (D) In the final step we construct a bijection between $\spec (n)$ and
$\mbox{tab}(n)$ (= set of all standard Young tableaux on the letters
$\{1,2,\ldots ,n\}$) sending
weights in $\spec (n)$ to content vectors of standard Young
tableaux showing, in particular, that the weights are integral.  
Weights in $\spec (n)$ that are related by $\sim$ 
go to standard Young tableaux of the
same shape. We shall not use this step. 

Let $V$ be a $S_n$-module, not necessarily multiplicity free. For $\alpha =
(a_1,a_2,\ldots ,a_n) \in \spec (n)$ define the {\em weight space}
$$V(\alpha) = \{v\in V\;:\;X_i(v)=a_iv,\;\;i=1,\ldots ,n\}.$$
Note that, if $V$ is multiplicity free, then it follows from items (B), (C)
above that every weight space is either
zero or one dimensional.

For $\lambda \in S_n^{\wedge}$, let $V(\lambda)\subseteq V$ denote the isotypical
component of $L^{\lambda}$ and let $\spec (n,\lambda)$ denote the set of all
weights $\alpha \in \spec (n)$ with $v_{\alpha}\in L^{\lambda}$. Basic
properties of the weight spaces are given below.
\bl \label{wsl}
Let $V,W$ be $S_n$-modules. We have

\noi (i) For $\lambda \in S_n^{\wedge}$, $V(\lambda)=\oplus_{\alpha\in \spec (n,\lambda)}V(\alpha)$
is an orthogonal
decomposition of $V(\lambda)$ under any
$S_n$-invariant inner product on $V$.

\noi (ii) $V = \oplus_{\alpha\in \spec (n)}V(\alpha)$ 
is an orthogonal
decomposition of $V$ under any
$S_n$-invariant inner product on $V$.

\noi (iii) Let $\lambda \in S_n^{\wedge}$, and $\alpha, \beta
\in \spec (n,\lambda)$. There is a canonical linear isomorphism
$f_{\alpha,\beta} : V(\alpha) \rar V(\beta)$, unique upto scalars,
satisfying the following property: for $v\in V(\alpha)$, the subspace
$\C [S_n]v \cap V(\beta)$ is
generated by $f_{\alpha,\beta}(v)$.

\noi (iv) Let $\lambda \in S_n^{\wedge}$, $\alpha \in \spec(n,\lambda)$, and
$f: V(\lambda) \rar W(\lambda)$ a $S_n$-linear map. Then 
$f(V(\alpha))\subseteq W(\alpha)$ and  
any linear map $g:V(\alpha)\rar W(\alpha)$ has a unique $S_n$-linear
extension $g:V(\lambda)\rar W(\lambda)$.
\el
\pf (i) Fix a $S_n$-invariant inner product on $V$ and let
$V(\lambda)=W_1\oplus \cdots \oplus W_t$ be an orthogonal decomposition of
$V(\lambda)$ into irreducible submodules, with all $W_i$ isomorphic to
$L^{\lambda}$. For $\alpha \in \spec(n,\lambda)$ it follows from (B), (C)
above that $V(\alpha)\cap W_i$ is one dimensional for all $i$ and that the
span of these one dimensional subspaces is $V(\alpha)$. Since the GZ-bases
of $W_i$ are orthogonal the result follows.

(ii) Follows from part (i), since the decomposition of $V$ into isotypical
components is orthogonal.

(iii) Let $v\in V(\alpha)$. Consider the irreducible submodule $\C[S_n]v$ of
$V$ generated by $v$. By item (B) above the subspace $\C[S_n]v \cap
V(\beta)$ is one dimensional, say generated by $u$. There is an element
$a\in \C[S_n]$ with $av=u$. The linear map 
\beq \label{gt}
V(\alpha)\rar V(\beta),\;x\mapsto ax
\eeq
has the required properties. It is also clear that
such a map is unique upto scalars.

(iv) The first part of the assertion is clear. For the second part note that
(\ref{gt}) implies that $g$ has atmost one $S_n$-linear extension. Dimension
considerations now show that $g$ has exactly one $S_n$-linear extension.
$\Box$

Consider the poset $M(n,k)$, on which the symmetric group $S_n$ acts by
substitution. Set $V=V(M(n,k))$ and $V_i = V(M(n,k)_i)$, $0\leq i \leq kn$.
Note that each $V_i$ is a $S_n$-submodule of $V$ and, for $\alpha \in \spec
(n)$, we have $V(\alpha)=V_0(\alpha) \oplus \cdots \oplus 
V_{kn}(\alpha)$. Since the up operator $U$ is $S_n$-linear it follows 
from Lemma \ref{wsl}(iv) that
\beq
\label{ws2}
U(V_i(\alpha)) & \subseteq &V_{i+1}(\alpha),\;\;0\leq i < kn,\;\alpha \in
\spec(n).
\eeq
Let $C=(v_i,v_{i+1},\ldots ,v_{kn-i})$, with $r(v_l)=l$ for all $l$, be a
symmetric Jordan chain in $V$. We say that $C$ is a {\em symmetric
Gelfand-Tsetlin chain} if each $v_l$ is a simultaneous eigenvector for the
action of $X_1,X_2,\ldots ,X_n$. It follows that $v_l\in V_l(\alpha)$,
$i\leq l \leq kn-i$, for some $\alpha \in \spec(n)$. Since an SJB of $V$
exists, it follows from Lemma \ref{wsl}(ii) and (\ref{ws2}) that there
exists a SJB of $V$ consisting of symmetric Gelfand-Tsetlin chains. However,
in such a SJB the Gelfand-Tsetlin chains belonging to different weight
spaces (of the same isotypical component) need not be related to each other.
Note that, for $\alpha \sim \beta,\;\alpha,\beta \in \spec(n)$, if $C$ is a
Gelfand-Tsetlin chain in $V(\alpha)$ then $f_{\alpha,\beta}(C)$ is a 
Gelfand-Tsetlin chain in $V(\beta)$ (from (\ref{gt}) and the fact that $U$ is
$S_n$-linear). We now add this condition to the definition.

A {\em symmetric Gelfand-Tsetlin basis} (SGZB) of $V$ is a SJB $B$ of $V$
satisfying the following conditions:

\noi (a) Each vector in $B$ is a simultaneous eigenvector for the action of
$X_1,\ldots ,X_n$.

\noi (b) For $\alpha \sim \beta,\; \alpha, \beta \in \spec(n)$, if $v\in B \cap V(\alpha)$, then
some multiple of $f_{\alpha,\beta}(v)$ is also in $B$.

Clearly, a SGZB of $V$ exists (for each isotypical component choose a SJB of
any one weightspace and transfer it to the other weightspaces via
$f_{\alpha,\beta}$). In the special case of Boolean algebras this definition
of SGZB coincides with the one given in the introduction (see proof of
Theorem \ref{mt3} below).

At this point two natural questions arise:

(i) Is it possible to characterize in some way the SJB produced by the
linear BTK algorithm? 

(ii) Is there an explicit construction (inductive or direct) of a SGZB of
$V(M(n,k))$?

Theorem \ref{mt3} answers both these questions in the special case of
Boolean algebras. We now prove this result. 

\bl \label{jl}
For $0\leq i \leq n$, $V(B(n)_i)$ is a multiplicity free $S_n$-module
with $\mbox{min}\{i,n-i\}+1$ irreducible summands.
\el

For the proof of this lemma see Theorem 29.13 in {\bf\cite{jl}}. To actually
identify the irreducibles (as corresponding to two part partitions) 
see Example 7.18.8 in {\bf\cite{s1}}. To identify
the irreducibles, along with their multiplicity (given by the number of
semistandard Young tableaux), in the $S_n$-module
$V(M(n,k)_i)$ see Exercise 7.75 in {\bf\cite{s1}}.

\pf {\bf of Theorem \ref{mt3}} We shall show inductively that each element
of $O(n)$ is a simultaneous eigenvector of $X_1,\ldots ,X_n$, the case $n=1$
being clear. By Lemma \ref{jl}, condition (b) in the defnition of SGZB is
then automatically satisfied. It also follows from Lemma \ref{jl} that in
the case of $V(B(n))$ the definition of SGZB given in the introduction
agrees with the definition given above.

Assume that each element of $O(n)$ is an eigenvector for the action of
$X_1,\ldots ,X_n$. Note that if $v\in V(B(n)_k)$ is an eigenvector for
$X_i$, for some $1\leq  i \leq n$,
then $\ol{v}\in V(B(n+1)_{k+1})$ is also an eigenvector for $X_i$ with the
same eigenvalue. Thus it follows from 
(\ref{ba1}, \ref{ba2}, \ref{ba3}, \ref{ba4}, \ref{ba5}) 
that each element of $O(n+1)$ is an eigenvector for $X_1,\ldots ,X_n$. It
remains to show that each element of $O(n+1)$ is an eigenvector for
$X_{n+1}$.

For $0\leq i \leq \frac{n+1}{2}$ and $0\leq k \leq i$
 define a subset $R(k,i)\subseteq O(n+1)$ consisting of all $v\in O(n+1)$
satisfying: $r(v)=i$ and the symmetric Jordan chain in which $v$ lies starts
at rank $k$ and ends at rank $n+1-k$. Put $W(k,i)= \mbox{Span }R(k,i)$. Clearly
$V(B(n+1)_i)=W(0,i) \oplus W(1,i) \oplus \cdots \oplus W(i,i)$. 
We claim that each $W(k,i)$ is a $S_{n+1}$-submodule. We prove this by
induction on $i$, the case $W(0,0)$ being clear. Assume inductively that
$W(0,i-1),\ldots ,W(i-1,i-1)$ are submodules, where $i\leq \frac{n+1}{2}$.
Since $U$ is $S_{n+1}$-linear, $U(W(j,i-1))=W(j,i)$, $0\leq j \leq i-1$ are
submodules. Now consider $W(i,i)$.
Let $u\in W(i,i)$ and $\pi\in S_{n+1}$. Since $U$ is
$S_{n+1}$-linear we have $U^{n+2-2i}(\pi u)=\pi U^{n+2-2i}(u)=0$. 
It follows that
$\pi u \in W(i,i)$. 

We now have from Lemma \ref{jl} that, for $0\leq i \leq \frac{n+1}{2}$,
 $W(0,i),\ldots ,W(i,i)$ are mutually
nonisomorphic irreducibles. Consider the $S_{n+1}$-linear map $f:V(B(n+1)_i)\rar
V(B(n+1)_i)$ given by $f(v) = av$, where
$$ a = \mbox{ sum of all transpositions in $S_{n+1}$ }= X_1+\cdots
+X_{n+1}.$$
It follows by Schur's lemma that there exist scalars $\lambda_0 ,\ldots
,\lambda_i$ such that $f(u)=\lambda_k u$, for $u\in W(k,i)$. Thus each element
of $R(k,i)$ is an eigenvector for $X_1+\cdots +X_{n+1}$ (and also for
$X_1,\ldots ,X_n$). It follows that each element of $R(k,i)$ is an
eigenvector for $X_{n+1}$. 

The paragraph above has shown that the bottom element of each symmetric
Jordan chain in $O(n+1)$ 
is a simultaneous  eigenvector for $X_1, \ldots ,X_{n+1}$. 
It now follows
from Lemma \ref{wsl}(iv) that each element of $O(n+1)$
is a simultaneous  eigenvector for $X_1, \ldots ,X_{n+1}$. 
That completes the proof. $\Box$

\begin{center} {\bf \Large{Acknowledgement}}
\end{center}
I am grateful to Sivaramakrishnan Sivasubramanian for several helpful
discussions and to Professor Alexander Schrijver for an encouraging e-mail.

\end{document}